\documentclass[10pt]{article}

\usepackage{amssymb, amsmath}

\newtheorem{Thm}{Theorem}
\newtheorem{Lem}[Thm]{Lemma}
\newtheorem{Cor}[Thm]{Corollary}
\newtheorem{Prop}[Thm]{Proposition}

\newtheorem{Rmk}{Remark}
\def\qed{\kern 6pt\hbox{\vrule\vbox to 6pt{\hrule width
                            6pt\vfil\hrule}\vrule}}

\newcommand{\N}{\mathbb{N}}

\newcommand{\R}{\mathbb{R}}
\newcommand{\C}{\mathbb{C}}

\def\ra{\rightarrow}

\def\eps{\varepsilon}

\def\L{\mathcal L}

\def\f{\varphi}

\title{Some functional forms of Blaschke-Santal\'o inequality}
\author{{ M. Fradelizi and
M. Meyer}}

\begin{document}

\maketitle

\begin{abstract}
We establish new functional versions of the Blaschke-Santal\'o inequality 
on the volume product of a convex body  which 
generalize to the non-symmetric setting an inequality of
K.~Ball \cite{BallPhD} and we give a simple proof of the case of equality. 
As a corollary, we get some inequalities  for 
$\log$-concave functions and Legendre transforms which extend the recent result of Artstein, 
Klartag and Milman \cite{AKM}, with its equality case. 
\end{abstract}

\vspace{3cm}
\noindent
Universit\'e de Marne la Vall\'ee,\\
Laboratoire d'Analyse et de Math\'ematiques Appliqu\'ees (UMR 8050)\\
Cit\'e Descartes - 5, Bd Descartes\\
Champs-sur-Marne\\
77454 Marne la Vall\'ee Cedex 2, France\\

\vspace{5mm}\noindent
Email: Matthieu.Fradelizi@univ-mlv.fr,\\
\hspace*{11mm}Mathieu.Meyer@univ-mlv.fr\\

\vspace{1mm}\noindent
Fax: 33 1 60 95 75 45

\newpage

\section{Introduction}
\noindent
For a Borel subset $K$ of $\R^n$ and a point $z\in \R^n$, the polar body
$K^{*z}$ of  $K$  with respect to $z$ is the convex set defined by:
$$K^{*z}=\{y\in \R^n; \langle y-z,x-z\rangle \le 1\hbox{ for every $x\in K$}\}.$$
Here $\R^n$ is endowed with the canonical scalar
product $\langle\ , \ \rangle$ and the associated  Euclidean norm $|\cdot|$. 
For $z=0$, we simply write $K^\circ$ instead of $K^{*0}$. 
Denote by $|A|$  the Lebesgue measure of a Borel subset $A$ of
$\R^n$. The {\it Santal\'o point} $s(K)$ of $K$ is a point for which
$$|K^{*s(K)}|=\min_{z}|K^{*z}|.$$ 
If $K$ is bounded and not contained in a hyperplane, its Santal\'o point $z$ is
characterized by the property that it is the center of
mass of $K^{*z}$. The inequality of Blaschke-Santal\'o (Blaschke
\cite{Blaschke},
Santal\'o \cite{Santalo}) states that
$$|K|\cdot |K^{*s(K)}| \le v_{n}^2:=|B_{2}^n|^2\ ,$$
where $B_{2}^n=\{x\in \R^n; |x|\le 1\}$ is the Euclidean ball.\\

We shall prove here new functional versions of the  Blaschke-Santal\'o inequality and 
give applications which extend the theorem  of Ball \cite{BallPhD}  as well as the
recent result of Artstein, Klartag and Milman \cite{AKM}. 
Notice that Lutwak and Zhang \cite{Lutwak1} and Lutwak,  Yang and Zhang \cite{Lutwak2}
gave other very different functional forms of the  Blaschke-Santal\'o inequality and  
recently Klartag and Milman \cite{KlartagMilman}, 
Klartag \cite{Klartag} and Colesanti \cite{Colesanti}  
also established functional forms of some other geometric inequalities. \\

The first main result of this paper generalizes 
with a new proof an inequality of K.~Ball \cite{BallPhD}; it treats the case of "centered" functions:

      \vskip 1mm \noindent
{\bf Proposition}  {\em Let $\rho:\R_{+}\to\R_{+}$ and $f_1, f_2:\R^n\to\R_+$ be
measurable functions such that
$$f_1(x)f_2(y)\le \rho^{2}(\langle x,y\rangle ) \hbox{ for every  $x,y\in
        \R^n$ satisfying $\langle x,y\rangle >0$}.$$
If the star shaped set $K_{1}=\{ x\in \R^n;  \int_{0}^{+\infty} r^{n-1} f_1(rx)dr\ge
        1\}$ is centrally symmetric (which holds if
        $f_{1}$ is even), or is a convex body with center of
        mass at the origin, then
         $$\int_{\R^n} f_1(x)dx\int_{\R^n} f_2(y)dy\le
        \left( \int_{\R^n} \rho({|x|^{2}})dx\right)^2.$$
 }
\vskip 1mm \noindent
The idea is to attach bodies $K_1$ and $K_2$ to the functions $f_1$ and $f_2$. From the duality relation on the 
$f_j$'s, we deduce, using the Pr\'ekopa-Leindler inequality for the geometric mean, that the 
sets $K_j$'s satisfy  the inclusion $K_2\subset c_n(\rho)K_1^\circ$, for some constant $c_n(\rho)$. 
Then the result follows from the Blaschke-Santal\'o inequality for sets. \\

As an application of this proposition, we treat the case of "non centered" functions:

\vskip 1mm \noindent
{\bf Theorem}  {\em Let $\rho:\R_{+}\to\R_{+}$ be measurable and $f:\R^n\to\R_+$ be
     a log-concave function such that $0<\int f <+\infty$. Then there
     exists $z\in \R^n$ with the following property: for any measurable function $ g: \R^n \mapsto \R_{+}$
     satisfying
     $$f(x)g(y)\le\rho^2\left(\langle x-z, y-z\rangle\right)\
$$
for every  $x,y\in \R^n$ with $ \langle x-z, y-z\rangle>0$,  one has
         $$\int_{\R^n} f(x)dx\int_{\R^n} g(y)dy\le
        \left( \int_{\R^n} \rho({|x|^{2}})dx\right)^2.$$
       }
\vskip 1mm \noindent
In the proof, we attach, for every $z\in\R^n$, the convex body 
$$K_{z}=\left\{ x\in \R^n;\ \int_{0}^{+\infty} f(z+rx) r^{n-1} dr \ge
1\right\}$$
and  show that there exists $z_{0}\in
\R^n$ such that the center of mass of $K_{z_{0}}$ is at the origin. Then the result
follows from  the preceding proposition. The existence of such a  $z_0$ is proved using
 Brouwer's fixed point theorem.\\

The main consequence of this  theorem is the following generalization of 
the results of Artstein, Klartag and Milman \cite{AKM} 
(who  considered only the cases $\rho(t)=e^{-t}$ and $\rho(t)=(1-t)_+^m$) 
for the Legendre transform $\L_{z}\phi$ of a convex function $\phi$.

\vskip 1mm\noindent
{\bf Theorem}  {\em Let $\rho:\R_{+}\to\R_{+}$ be a log-concave
non-increasing function and let $\phi$ be a convex
function such that
$0<\int_{\R^n}\rho\left(\phi(x)\right)dx<+\infty$ . Then for some $z\in
\R^n$,
one has $$\int_{\R^n}\rho\left(\phi(x)\right)dx
\int_{\R^n}\rho\left(\L_{z}\phi(y)\right)dy\le
\left(\int_{\R^n}\rho\left(\frac{|x|^{2}}{2}\right)dx\right)^2.$$
}
\vskip 1mm \noindent

In all these functional forms of Blaschke-Santal\'o inequality, we determine the equality cases and 
establish some geometric corollaries. In particular we investigate the following question:

What are the Borel measures $\mu$ on $\R^n$ and the sets $K$ in $\R^n$ which 
satisfy a Blaschke-Santal\'o type inequality
$$\mu(K)\cdot \mu(K^\circ)\le \mu(B_2^n)^2\ ?$$
Cordero-Erausquin (\cite{Cordero}) proved such an inequality 
in $\C^n$ for plurisubharmonic measures and $\C$-symmetric 
pseudo-convex sets, using complex interpolation. 
He also remarked that it holds
for the Gaussian measure in $\R^n$ and asked whether it  still holds
for any symmetric $\log$-concave measures $\mu$ and  any symmetric convex body $K$  in $\R^n$. 
Klartag also established this inequality for a special class of measures in \cite {Klartag}.
As corollaries of our functional inequalities, we get that this inequality holds:\\
- for any unconditional  $\log$-concave measure $\mu$ and unconditional measurable set $K$\\
- for any rotation invariant $\log$-concave measure $\mu$ and any centrally symmetric measurable set $K$.\\
And we determine the equality cases.

\vskip 2mm\noindent
The paper is organized in the following way. In section 2, we  treat the case of
unconditional functions and sets, where one can apply a multiplicative
version of the Pr\'ekopa-Leindler inequality. In Section 3, we prove the proposition stated above concerning the
case of "centered" functions. Section 4 is
devoted to the proof of our  theorem on general (not centered) functions. 
In Section 5, we prove the consequences for Legendre
transforms of convex functions.
\vskip 2mm\noindent
It should be observed that the main difficulty when working with
Santal\'o type inequalities for non-symmetric bodies or functions is to find a good
center. If $G(K)$ is the center of mass of $K$ ($G(K)=\int_{K}xdx/|K|$), 
one has as well
$$|K|\cdot |K^{*G(K)}| \le v_{n}^2, $$
because  Blaschke-Santal\'o  inequality can be applied to
$K^{*G(K)}$. But if $K$ is centrally symmetric, the situation is simpler:
$\min_{z}|K^{*z}|$ is reached at $0$, and then $|K|\cdot  |K^\circ| \le
|B_{2}^n|^2$.
We shall also make use of the equality case in Blaschke-Santal\'o inequality:  
there is equality if and only if $K$ is an ellipsoid. 
At the end of the paper, we give a new and elementary proof of this result.

\vskip 3mm \noindent
\section  {An inequality for unconditional functions}


       \vskip 2mm \noindent
       We say that a function $\f: \R^n\mapsto \R$ is {\it unconditional} if
       $$\f(\eps_1 x_1, \dots, \eps_n x_n )=\f( x_1, \dots, x_n)$$ 
       for every $(\eps_1,\dots,\eps_n)\in \{-1,1\}^n$ and every $(x_1,\dots,x_n)\in\R^n$.
In the same way, a subset $K$ in $\R^n$ is {\it unconditional }
if its characteristic function $\chi_{K}$ is
unconditional. Observe that an unconditional convex function
$W:\R^n\mapsto \R$ is
minimal at $0$ and is moreover   {\it increasing}, in the sense that
$W(x)\le W(y)$ whenever
$x=(x_{1},\dots,x_{n})$ and $y=(y_{1},\dots,y_{n})$ satisfy $|x_{i}|\le
|y_{i}|$, $1\le i\le n $.
\vskip 1mm \noindent
In particular, if $W$
is unconditional and convex, one has
$$
W(\sqrt{x_1 y_1},\dots, \sqrt{x_n y_n})\le W\left(\frac{x+y}{
2}\right)\le \frac{W(x)+W(y)}{2}\  , $$
for all $x=(x_{1},\dots,x_{n}),y=(y_{1},\dots,y_{n})\in\R_+^n$

\vskip 3mm \noindent
The next proposition is a  form of Pr\'ekopa-Leindler inequality for the
geometric mean due to Borell (\cite{Borell}), Ball (\cite{Ball}),
Uhrin (\cite{Uhrin}).
This result is well known and follows  from the usual
Pr\'ekopa-Leindler
inequality. We prove it here for
the convenience of the reader.
As we shall see in the corollary, this proposition gives
a first functional form of Blaschke-Santal\'o inequality.

\begin{Prop}\label{Prekopageom}{\bf (Pr\'ekopa-Leindler inequality for
the geometric mean)}\\
Let $f_1$, $f_2$, $f_3:\ \R^n\to\R_{+}$  be unconditional measurable
functions such that
$$ f_1(x_1,\dots,x_n)
f_2(y_1,\dots,y_n)\le f_3(\sqrt{x_1 y_1},\dots, \sqrt{x_n y_n})^2$$
for every $(x_1,\dots,x_n)$ and $(y_1,\dots,y_n)\in \R_+^n$. Then
$$\int_{\R^n} f_1(x)dx\int_{\R^n} f_2(y)dy\le \left(\int_{\R^n} f_3(z)dz\right)^2$$
with equality if and only if there exists a continuous function $\tilde{f_3}:\R_{+}\to\R_{+}$ such that the following two conditions hold:
  \vskip 1mm \noindent
{\bf a.} $f_3=\tilde{f_3}$ a.e. and $\tilde{f_3}(x_1,\dots,x_n)\tilde{f_3}(y_1,\dots,y_n)\le \tilde{f_3}(\sqrt{x_1 y_1},\dots, \sqrt{x_n y_n})^2$ 

        \vskip 1mm \noindent
{\bf b.} for some $c_1,\dots, c_n>0$ and $d>0$, one has 
$$f_1(x_1, \dots,x_n)=d\tilde{f_3}(c_1 x_1, \dots,c_nx_n)\hbox{ and }  f_2(x) 
=\frac{1}{d}
\tilde{f_3}\left( \frac{x_1}
{c_1}, \dots,\frac{x_n}{c_n}\right) \quad  a.e.$$
      
\end{Prop}
        \vskip 5mm \noindent
        {\bf Proof:} Since the $f_{j}$ are unconditional,
one has  $ \int _{\R^n}f_j= 2^n \int_{\R_+^n} f_{j}$, \ $j=1,2,3$. 
For $(t_1,\dots, t_n)\in\R^n$, we define
        $$g_j(t_1,\dots, t_n)=f_j(e^{t_1}, \dots, e^{t_n})\,e^{\sum_{i=1}^n t_i}\ .$$ 
We get
        $$\int_{\R_+^n} f_j =\int_{\R^n} g_j$$
and for every $s,t\in \R^n$,
        $$g_1(s) g_2(t)\le g_3\left(\frac{s+t}{2}\right)^2.$$
Hence the result follows from  Pr\'ekopa-Leindler inequality. For the equality
case, see \cite{Dubuc}.
\hfill\qed

       \vskip 5mm \noindent
As a corollary, we get the following generalized form of
Blaschke-Santal\'o inequality for
unconditional sets, together with its case of equality.

\begin{Cor}\label{corinc}
 Let $W: \R^n\to\R\cup\{+\infty\}$ be an unconditional convex function 
 and  let $\mu$ be the Borel measure on $\R^n$ with density $e^{-W(x)}$ with
respect to the Lebesgue measure. Then one has
$$
\mu(K)\mu(K^\circ)\le\mu(B_2^n)^2,
$$
for every unconditional measurable set $K\subset\R^n$.

\noindent
If moreover the support of $\mu$ is $\R^n$, there is equality if and
only if there exists a diagonal matrix $T$, with diagonal entries $(t_1,...,t_n)\in\R_+^n$
such that:

- $K=T(B_2^n)$

-   $W(x)=W(Px)$,  for every $ x\in K\cup K^\circ\cup B_2^n$, where $P$  is the orthogonal projection on  the subspace spanned by the $(e_i)_{i\in I}$ and $I=\{i ;\  1\le i\le n,\ t_i=1\}$.
\end{Cor}

\vskip 3mm \noindent
{\bf Proof:}
\vskip 1mm \noindent
{\bf A. The inequality.}
\vskip 1mm
\noindent
We apply Proposition \ref{Prekopageom}  to
$$
f_1(x)=e^{-W(x)}\chi_K(x),\ f_2(x)=e^{-W(x)}\chi_{K^\circ}(x),
f_3(x)=e^{-W(x)}\chi_{B_2^n}(x)\ .
$$
The hypotheses are satisfied since for all $x=(x_1,\dots,x_n),
y=(y_1,\dots,y_n)\in \R_+^n$, one has
$$
\chi_K(x)\chi_{K^\circ}(y)\le\chi_{B_2^n}(\sqrt{x_1 y_1},\dots,
\sqrt{x_n y_n})\
$$
and
\begin{equation}\label{Wconv}
W(\sqrt{x_1 y_1},\dots, \sqrt{x_n y_n})\le W\left(\frac{x+y}{2}\right)\le
\frac{W(x)+W(y)}{2}\  .
\end{equation}
\noindent
as explained at the beginning of this section.
This gives the inequality.
\vskip 2mm
\noindent
{\bf B. The case of equality.}
\vskip 1mm
\noindent
Assume that the support of $\mu$ is $\R^n$ (hence  $W(x)< +\infty$, for every $x\in\R^n$) and that there is equality
in the preceding inequality. From the
equality case in Proposition \ref{Prekopageom},  there exists
$t_1,\dots, t_n>0$ and $d>0$, such that if we
denote  by $T$ the diagonal
matrix with diagonal entries  $(t_1,\dots, t_n)$, then
$$
e^{-W(x)}\chi_K(x)=de^{-W\left(Tx\right)}\chi_{B_2^n}\left(Tx\right)
$$
and
$$
\ e^{-W(x)}\chi_{K^\circ}(x)=\frac{1}{
d}e^{-W\left(T^{-1}x\right)}\chi_{B_2^n}\left(T^{-1}x\right) .
$$
We get $K=T^{-1}(B_2^n)$ and $K^\circ=T(B_2^n)$.
Taking $x=0$ gives $d=1$ so that
$$
W(x)=W\left(Tx\right)=W\left(T^{-1}x\right)\ \hbox{for every $x\in B_2^n$ .}
$$
Let $S=\frac{T+T^{-1}}{ 2}$ be the diagonal matrix with  diagonal entries
$s_i=\frac{1}{ 2}\left(t_i +\frac{1}{ t_i}\right)$, $1\le i\le n$.
One has  $s_i>1$ for all $i\notin
I:=\{j\ ;\ t_j=1\}$ hence $\lim_{k\rightarrow +\infty}S^{-k}(x)=Px$, for all $x\in\R^n$.
Using the inequalities (\ref{Wconv}) for $Tx$ and $T^{-1}x$, we get
$$
W(x)\le W\left(\frac{Tx+T^{-1}x}{2}\right)\le\frac{W(Tx)+W(T^{-1}x)}{ 
2}=W(x)\  .
$$
 Hence $W(Sx)=W(x)$ for every $x\in B_2^n$.
The result follows from the continuity of $W$.
\hfill\qed

\vskip 2mm \noindent
{\bf Remarks:}\\
{\bf 1) } Actually the proof shows that 
the inequality of Corollary \ref{corinc} still holds true when the 
hypothesis that $W$ is convex is replaced with the weaker hypothesis that
$$ (t_1,\dots,t_n)\mapsto W(e^{t_1},\dots,e^{t_n})$$
is convex on $\R^n$.\\
{\bf 2) }The Pr\'ekopa-Leindler inequality for the geometric mean was also used in
\cite{CFM} to prove that if $K$ is an unconditional convex body and
$\mu$ has an unconditional log-concave density with
respect to the Lebesgue measure, then $t\mapsto\mu(e^t K)$ is a log-concave function.


\section{The Blaschke Santal\'o inequality for centered functions.}

In the next result, we generalize with a new proof an inequality
obtained by K. Ball \cite{BallPhD}
in the special case of even functions, and we characterize the case 
of equality.

      \begin{Prop}\label{even} Let $\rho:\R_{+}\to\R_{+}$ and $f_1, f_2:\R^n\to\R_+$ be
measurable functions such that
$$f_1(x)f_2(y)\le \rho^{2}(\langle x,y\rangle ) \hbox{ for every  $x,y\in
        \R^n$ satisfying $\langle x,y\rangle >0$}\ .$$
If the star shaped set $K_{1}=\{ x\in \R^n;  \int_{0}^{+\infty} r^{n-1} f_1(rx)dr\ge
        1\}$ is centrally symmetric (which holds if
        $f_{1}$ is even), or if $K_1$ is a convex body with center of
        mass at the origin, then
         $$\int_{\R^n} f_1(x)dx\int_{\R^n} f_2(y)dy\le
        \left( \int_{\R^n} \rho({|x|^{2}})dx\right)^2$$
with  equality if and only if  for some continuous function $\tilde{\rho}:\R_{+}\to\R_{+}$ one has 
\vskip 1mm\noindent
{\bf a.} $\rho=\tilde{\rho}$ a.e., $\sqrt{\tilde\rho(s)\tilde\rho(t)}\le\tilde\rho(\sqrt{st})$ for every $s,t\geq 0$ 
and if $n\geq 2$, $\tilde\rho(0)>0$ or $\tilde\rho$ is the null function.
\vskip 1mm \noindent
{\bf b.}  For some positive definite $[n\times n]$ matrix $T$ and for some $d>0$, one has 
$$f_1(x) = d\tilde\rho( |Tx|^2)\hbox{ and }  f_2(x) =\frac{1}{d}\tilde\rho(|T^{-1}x|^2)\quad a.e.$$
 \end{Prop}

\noindent
         {\bf Proof:}
\vskip 1mm \noindent
{\bf A. The inequality.}
\vskip 1mm
\noindent
Let $x_1, x_2\in \R^n$ satisfying $\langle x_1,x_2\rangle >0$. We
define $g_j: \R_+\ra \R_+$ by
$$g_j(s)=
       s^{n-1} f_j(sx_j),\  j=1,2\ \hbox{ and }g_3(u)= u^{n-1} \rho 
(u^{2}\langle
       x_1,x_2\rangle).$$
       Then by hypothesis, one has
        $g_1(s)g_2(t)\le (st)^{n-1} \rho^{2}(st \langle x_1,x_2\rangle)=
        g_3^2(\sqrt{st}) \ .$
        It follows from  Proposition \ref{Prekopageom} ($n=1$) that
      $$
      \int_{\R_+} s^{n-1} f_1(sx_1)ds\int_{\R_+}t^{n-1} f_2(tx_2)dt
\le \left(\int_{\R_+} u^{n-1}\rho\left((u^{2}\langle x_1,x_2\rangle\right)
du\right)^2$$
$$=
\frac{1}{\langle
x_1,x_2\rangle^{n}}\left(\int_{\R_+}r^{n-1}\rho(r^{2}))
dr\right)^2= \frac{c_{n}(\rho)^n}{\langle x_1,x_2\rangle^{n}} .
$$
where $c_n(\rho): =\left(\int_{\R_+} r^{n-1}\rho(r^{2}) dr\right)^{\frac{2}{ n}}$.
For $j=1,2$, we define
        $$K_{j} =\{x\in \R^n; \int_{\R_+} r^{n-1} f_j(rx)dr\ge 1\}\ .$$
The sets $K_1$ and $K_2$ are starshaped with respect to the origin.
Denote their gauge by $\|\cdot\|_{K_j}$, $j=1,2$. One has
$$
\|x\|_{K_j}=\inf\{ \lambda >0; \ x\in \lambda K_j\} =
\left(\int_{\R_+} r^{n-1} f_j(rx)dr\right)^{-{\frac{1}{ n}}}
\hbox {for all $x\in\R^n$}\ .$$
The preceding inequality may be read as follows: for every $x_1, x_2\in \R^n$
such that $\langle x_1,x_2\rangle >0$, one has
\begin{equation}\label{polar}
\langle x_1,x_2\rangle\le c_n(\rho)\|x_1\|_{K_1}\|x_2\|_{K_2}\ .
\end{equation}
This means that
$$K_2\subset c_n(\rho)K_1^\circ\ .$$

\medskip\noindent
Under our hypotheses, either $K_{1}$ is centrally symmetric, so its
closed convex hull is also centrally symmetric and has its center of mass at
the origin, or $K_{1}$ is itself a convex body with center of
mass at the origin.  In both cases, the origin is actually the Santal\'o point of
$K_{1}^{\circ}$, and it follows from Blaschke-Santal\'o inequality that
$|K_{1}|\ |K_{1}^\circ|\le v_{n}^{2}$.
We get thus
$$|K_1|\ |K_2|\le c_n(\rho)^n |K_1|\ |K_1^\circ|\le c_n(\rho)^nv_n^2\ .$$
Integrating in polar coordinates for $j=1,2$, one has
$$
\int_{\R^n}f_j(x)dx=nv_n\int_{S^{n-1}}\int_{\R_+}s^{n-1}f_j(su)dsd\sigma(u)
=nv_n\int_{S^{n-1}}\frac{d\sigma(u)}{\|u\|_{K_j}^n}=n|K_j|\ ,
$$
where $\sigma$ denotes the rotation invariant probability on the unit sphere 
$S^{n-1}:=\{u\in\R^n\ ;\ |u|=1\}$. Thus
$$
\int_{\R^n} f_1(x)dx\int_{\R^n} f_2(y)dy=n^2|K_1||K_2|\le (nv_n)^2 c_n(\rho)^n=
\left( \int_{\R^n} \rho({|x|^{2}})dx\right)^2.
$$
{\bf B. The case of equality.}
\vskip 1mm
\noindent
Assume now that there is equality. By the case of equality of
Blaschke-Santal\'o inequality, $K_1$ is an
ellipsoid centered at the origin and $K_2=c_n(\rho)K_1^\circ$. We may and do assume that $K_1=B_2^n$. 
For every $x\in S^{n-1}$, one has
$\langle x,x\rangle =1=c_n(\rho)\|x\|_{K_1}\|x\|_{K_2}$,
which means that there is equality in (\ref{polar}) for $x_1=x_2=x$.
From the equality case  of Proposition \ref{Prekopageom} ($n=1$), it follows
 that there exists a continuous function $\tilde{\rho}:\R_{+}\to\R_{+}$ such that
 \vskip 1mm
\noindent
-  $\rho=\tilde{\rho}$ a.e., $\sqrt{\tilde\rho(s)\tilde\rho(t)}\le\tilde\rho(\sqrt{st})$ for every $s,t\geq 0$ 
 \vskip 1mm
\noindent
- for every $x\in S^{n-1}$, there exists $c=c(x)>0, d=d(x)>0$ such that
$$
g_1(s)=dg_3(cs)\hbox{ and }g_2(s)=\frac{1}{ d} g_3\left(\frac{s}{c}\right) 
\hbox{ for a.e. $s\ge 0$.}
$$
Let us prove that $c$ and $d$ are constant functions. Since
$$
1=|x|^{-n}=\|x\|_{K_1}^{-n}
=\int_{\R_+} g_1(s)ds=\frac{d(x)}{c(x)}\int_{\R_+}
g_3(u)du=(c_n(\rho))^{\frac{n}{2}}\frac{d(x)}{c(x)}\ ,
$$
we have $d(x)=\frac{c(x)}{ c_n(\rho)^{n/2}}$. Hence for a.e.  $s\ge 0$
$$f_1(sx)=\left(\frac{c(x)}{\sqrt{c_n(\rho)}}\right)^n\tilde\rho(c(x)^2s^2)\ ,\
f_2(sx)=\left(\frac{\sqrt{c_n(\rho)}}{c(x)}\right)^n\tilde\rho\left(\frac{s^2}{c(x)^2}\right)
$$
By the hypotheses, for every $x,y\in S^{n-1}$ 
satisfying $\langle x,y\rangle >0$ and  $s,t\geq 0$
$$
\left(\frac{c(x)}{c(y)}\right)^n\tilde\rho(c(x)^2s^2)\tilde\rho\left(\frac{t^2}{c(y)^2}\right)\le\tilde\rho^2(st\langle x,y\rangle)\ .
$$
If $\tilde\rho(0)\neq 0$, we take  $s=t=0$, simplify and get $c(x)\le c(y)$, for any $x,y\in S^{n-1}$. Therefore  $c$ is a constant function. \\
If $\tilde\rho(0)= 0$ and $n\ge 2$, we take $x, y\in S^{n-1}$ with $\langle x,y\rangle =0$ 
(this  is possible since $\tilde\rho$ is continuous), we get that  $\tilde\rho$ is the null function.
\hfill\qed
       \vskip 3mm \noindent
       {\bf Remarks:}
    \vskip 1mm \noindent
{\bf 1)} We did not follow here the  more natural proof
       given  by K.~Ball in the even case. For sake of completeness,
   we outline his  proof in the case where $\rho$ is non-increasing.  
   Setting for $t>0$, $i=1,2$, $p_{i}(t)=|\{f_{i}>t\}|$,
       one has $\int f_{i}=\int_{0}^{+\infty} p_{i}(t)dt$.
    The hypothesis on $f_{1}$ and $f_{2}$ gives that for every $s,t>0$,
    one has $\{f_{2}>t\}\subset \rho^{-1} (\sqrt{st})\{f_{1}>s\}^{\circ}$.
   Now, {\it the fact that $f_{1}$ is even} implies that its
    level sets are centrally symmetric and this allows to apply
    Blaschke-Santal\'o inequality
    to get  for all $s,t>0$,
$$p_{1}(s)p_{2}(t)\le \left(\rho^{-1} (\sqrt{st})\right)^n 
v_{n}^{2}\, ,$$
and the result follows from
Proposition~\ref{Prekopageom} applied in dimension $1$.
\vskip 1mm \noindent
{\bf 2)} The idea of attaching a convex set of the form of $K_1$ to a
$\log$-concave function $f_{1}$ to prove a functional inequality was
originally used by K.~Ball in \cite{Ball} and is also used by Klartag 
and Milman in \cite{KlartagMilman}.
\vskip 1mm \noindent
{\bf 3)} There are many ways to recover the usual Blaschke-Santal\'o inequality
for symmetric sets from Proposition~\ref{even}. As noticed by
K.~Ball in \cite{BallPhD}, the more natural is  to apply it
to $f_1=\chi_K$, $f_2=\chi_{K^\circ}$ and $\rho=\chi_{[0,1]}$.
But more generally,  we get the same result by
applying it to $f_1(x)=\rho(\|x\|_K^2)$, $f_2(y)=\rho(\|y\|_{K^\circ}^2)$
and any function $\rho$ such that
$t\mapsto\rho(e^t)$ is log-concave and non-increasing on $\R$.
This was noticed by Artstein, Klartag and Milman \cite{AKM} in the 
case when $\rho(t)=e^{-t}$.
\vskip 1mm \noindent
{\bf 4)} Let $K$ be a convex body whose center of mass is at the origin. If
we set $f_1=\chi_K$, $f_2=\chi_{K^\circ}$ and $\rho=\chi_{[0,1]}$,
we get $K_{1}=K/n^{1/n}$ so that center of mass of $K_{1}$ is at the origin.
Hence Proposition~\ref{even} also permits to recover 
the general Blaschke-Santal\'o inequality
for convex sets.

\vskip 3mm \noindent
As a corollary of Proposition~\ref{even}, let us prove a generalized form of
Blaschke-Santal\'o inequality for
symmetric sets and some class of rotation invariant measures.
This inequality is known for the Lebesgue measure and the Gaussian measure
(see \cite{Cordero}); and also for a special class of measures 
(see \cite{Klartag}).
It was asked in \cite{Cordero} whether it holds
for any symmetric $\log$-concave measure. We also give here a partial
answer:

\begin{Cor} Let $h:\R_+\to\R_+$
be a non-increasing function which satisfies that
      $t\mapsto h(e^t)$ is $\log$-concave on $\R$. Let $\mu$ be the rotation
invariant measure on $\R^n$,
   with density
      $h(|x|)$ with respect to the Lebesgue measure.
   Then, for every centrally symmetric measurable set
      $K\subset\R^n$, one has
$$
\mu(K)\mu(K^\circ)\le\mu(B_2^n)^2.
$$
If moreover, the support of $\mu$ is $\R^n$, there is equality if and only if

- either $K=B_2^n$

- or $K=T(B_2^n)$ for some positive definite matrix $T\neq I$ and
$h$ is constant on $[0,\max(\|T\|,\|T^{-1}\|)]$, where $\|T\|=\max_{|x|=1}|Tx|$.
\end{Cor}

\noindent
         {\bf Proof:}
         \vskip 1mm \noindent
{\bf A. The inequality.}
\vskip 1mm
\noindent
We apply Proposition~\ref{even} to
$$
f_1(x)=h(|x|)\chi_K(x),\ f_2(y)=h(|y|)\chi_{K^\circ}(y)\ {\rm and}\
\rho(t)=h(\sqrt{t})\chi_{[0,1]}(t).
$$
The hypotheses are satisfied since for all $x,y\in \R^n$ such
that $\langle x,y\rangle>0$, one has
$$
f_1(x)f_2(y)\le h^2\left(\sqrt{|x||y|}\right)\chi_{[0,1]}(\langle x,y\rangle)
\le h^2(\sqrt{\langle x,y\rangle})\chi_{[0,1]}(\langle x,y\rangle)=
   \rho^2\left({\langle x,y\rangle}\right)$$
and $f_1$ is even. We get thus
$$
\int f_{1}(x)dx\int f_{2}(y)dy=\mu(K)\mu(K^\circ)\le
\left( \int_{\R^n} \rho({|x|})dx\right)^2=\mu(B_2^n)^2.
$$
{\bf B. The case of equality.}
\vskip 1mm
\noindent
Assume that the support of $\mu$ is $\R^n$ (hence $h>0$) and that there is equality.
It follows from  Proposition~\ref{even} that for
some positive matrix $T$
and for some $d>0$, one has
$$f_1(x) = d\,\rho( |Tx|^{2})\hbox{ and
}f_2(y) =\frac{1}{d}\rho(|T^{-1}y|^{2})\hbox{ for all $x,y\in\R^n$\ .}$$
This gives
$$
h(|x|)\chi_K(x)=d\,h(|Tx|)\chi_{[0,1]}(|Tx|)\ $$
and $$
 h(|y|)\chi_{K^\circ}(y)=\frac{1}{d}h(|T^{-1}y|)\chi_{[0,1]}(|T^{-1}y|)\ .
 $$
Hence $K=T^{-1}(B_2^n)$, $K^\circ=T(B_2^n)$ and 
$h(|Tz|)=h(|z|)=h(|T^{-1}z|)$ for every $z\in B_2^n$.
If $K\neq B_2^n$, one has $\max(\|T\|,\|T^{-1}\|)>1$.
We may assume that $\|T\|>1$.
Let $z_0\in S^{n-1}$ satisfying $|Tz_0|=\|T\|$ and $\lambda\in [0,\|T\|]$.
Applying the previous equality to
$z={\lambda z_0/\|T\|}$, we get
$h(\lambda)=h(|Tz|)=h(|z|)=h(\lambda/\|T\|) .$
From the continuity of $h$,
$h(\lambda)=h(\lambda/\|T\|^n) =h(0).$
\hfill\qed


\section{The general case}

We are now in position to prove the following theorem.

       \vskip 3mm \noindent
\begin{Thm}\label{Main} Let $\rho: \R_{+}\to\R_{+}$ be measurable and
     $f: \R^n\to\R_+$ be
     a log-concave function such that $0<\int f <+\infty$. Then there
     exists $z\in \R^n$ such that for any measurable function $ g: \R^n \mapsto \R_{+}$
     satisfying
     $$f(x)g(y)\le\rho^2\left(\langle x-z, y-z\rangle\right)\ $$
for every  $x,y\in \R^n$ such that $ \langle x-z, y-z\rangle>0$,  one has
         $$\int_{\R^n} f(x)dx\int_{\R^n} g(y)dy\le
        \left( \int_{\R^n} \rho({|x|^{2}})dx\right)^2$$
with  equality if and only if the following two conditions hold:
\vskip 1mm\noindent
{\bf a.}  For some positive definite $[n\times n]$ matrix $T$, some $z\in\R^n$ and some $d>0$, 
$$f(x) = d\rho\left( |T(x-z)|^2\right)\hbox{ and }  g(x) =\frac{1}{d}\rho\left(|T^{-1}(x-z)|^2\right)\quad a.e.$$
\vskip 1mm \noindent
{\bf b.} $\sqrt{\rho(s)\rho(t)}\le\rho(\sqrt{st})$ a.e.
 \end{Thm}
 
  \vskip 2mm\noindent
   {\bf Proof:}
   \vskip 1mm\noindent
   For every $z\in \R^n$ let
$$K_{z}=\left\{ x\in \R^n;\ \int_{0}^{+\infty} f(z+rx) r^{n-1} dr \ge
1\right\}.$$
Since $f$ is $\log$-concave, it follows from Ball \cite{Ball} that for
every $z\in\R^n$,
the set $K_{z}$ is a convex body. If we can prove that there exists $z_{0}\in
\R^n$ such that the center of mass of $K_{z_{0}}$ is at the origin, we get the result
from  proposition \ref{even} applied to $f_1(x)=f(x+z_0)$ and $f_2(x)= g(x+z_0)$.
\hfill\qed
   \vskip 2mm \noindent
This will be done in the following two lemmas, using Brouwer's fixed point theorem.

   \vskip 2mm \noindent
\begin{Lem}\label{lemr} Let $n\ge 2$ and  $f:\R^n\to \R_{+}$ be a log-concave function
such that $0< \int f <+\infty$.
For $z,x\in \R^n$,
define $r_{z}(x)=\left(\int_{0}^{+\infty} f(z+rx) r^{n-1} dr\right)^{\frac{1}{n}}.$ 
One has then
\vskip 1mm\noindent
{\bf 1)} For all $\eps>0$ and $\alpha <1$, there exists $M>0$ such that
$r_{z}(u)\le \eps$
whenever
$u\in S^{n-1}$ and $z\in \R^n$ satisfy   $\langle u,z\rangle \ge -\alpha |z|$
and $|z|\ge M$.
\vskip 1mm\noindent
{\bf 2)} $r_{z}\left(-\frac{z}{|z|}\right)\ra +\infty$ when $|z|\ra +\infty$.
\end{Lem}

\vskip 2mm\noindent
{\bf Proof:}
\vskip 1mm\noindent
 From the hypotheses on $f$, it is easy to see that for some $a,b,c,d>0$,
one has
$$a\chi_{bB_{2}^n}(x)\le f(x)\le d\,e^{-c|x|}\hbox{ for every } x\in
\R^n\ .$$
\vskip 1mm\noindent

1)  If $u\in S^{n-1}$ and $z\in \R^n$ satisfy \
$-\langle u,z\rangle \le \alpha |z|\ , $ then for every $r\ge 0$
$$|z+ru|^2 \ge |z|^{2}-2\alpha |z| r +r^2
\ge (1-\alpha)(|z|^{2}+r^{2})\ge \frac{1-\alpha}{ 2} (|z|+r)^2.$$
It follows that
$$ r_{z}(u)^n \le  de^{ -c\sqrt{\frac{1-\alpha}{2} }|z|}
\int_{0}^{+\infty} r^{n-1} e^{-c\sqrt{\frac{1-\alpha}{2}} r}dr\ra 0\
\hbox{when}\ |z|\ra +\infty. 
$$
\vskip 1mm\noindent
2)  Let $u=-\frac{z}{ |z|}$. Then
$$     r_{z}(u)^n   =\int_{0}^{+\infty} r^{n-1}f\left((r-|z|)u\right) dr\ge
a\int_{0}^{+\infty}r^{n-1} \chi_{[-b,b]}(r-|z|)\, dr$$
Thus, for $|z|>b$, one has
$r_{z}(u)^n  \ge \frac{a}{ n}\left((|z|+b)^n
-(|z|-b)^n\right)\ra +\infty$ when  $|z|\ra +\infty$.\hfill\qed

\vskip 2mm\noindent
As we have already seen, for every $z\in\R^n$, the set $K_z$ is a convex body. 
Moreover notice that under our hypotheses, the origin is in the interior of $K_z$ and 
$r_z$ is the radial function of $K_z$ 
($r_z(u)=\max\{\lambda>0\ ;\ \lambda u\in K_z\}$, for every $u\in S^{n-1}$). Hence part 1) of the 
preceding lemma means that for all $\eps>0$ and $\alpha <1$, there exists $M>0$ such that
for every $|z|\ge M$,
$$\{x\in K_z\ ;\ \langle x,z\rangle\ge -\alpha |z|\}\subset \eps B_2^n\, .$$

\vskip 3mm\noindent
\begin{Lem} Let $f: \R^n\to \R_{+}$ be a log-concave function
such that $0<\int f <+\infty.$ For every $z\in \R^n$ let
$$K_{z}=\left\{ x\in \R^n;\ \int_{0}^{+\infty} f(z+rx) r^{n-1} dr \ge
1\right\}.$$
Then there exists  $z_{0}\in \R^n$ such that the convex body $K_{z_{0}}$ has its 
center of mass  at the origin.
\end{Lem}

\vskip 1mm\noindent
{\bf Proof:}
\vskip 1mm\noindent
Notice first that for $n=1$, the result is easy, one chooses the unique point $z_0\in\R$ such that
$$\int_{z_0}^{+\infty}f(r)dr=\int_{-\infty}^{z_0}f(r)dr,$$
then $K_{z_0}$ is a symmetric interval. We assume from now on that $n\ge 2$.
It is clear that $z\mapsto K_{z}$ is continuous for the
Hausdorff distance, so that if $G(z)$ is the centre of mass
of $K_{z}$, then $G:\R^n\mapsto \R^n$ is continuous.

\vskip 1mm\noindent
{\bf A.} We first show  that
$$\ |G(z)|\ra +\infty\ \hbox{ and }
\ \big\langle \frac{G(z)}{|G(z)|}, \frac{z}{ |z|}\big\rangle \ra -1\
\hbox{ when }|z|\ra +\infty\ .$$
\vskip 1mm\noindent
Let $h_{K_{z}}$ be the support function of
$K_{z}$ {\em i.e.} 
$$h_{K_{z}}(y)=\max_{x\in K_{z}}\langle x,y\rangle \hbox{ for every }y\in \R^n.$$
It is well known that one has, for all $u\in S^{n-1}$,
$$-h_{K_z}(-u) +\frac { h_{K_z}(u)+h_{K_z}(-u)}{n+1})\le
\langle G(z),u\rangle \le h_{K_z}(u)-\frac{ h_{K_z}(u)+h_{K_z}(-u)}{n+1}\ .$$
By part 1) of Lemma \ref{lemr} applied with $\alpha =0$, 
for every $\eps >0$, there exists $M>0$ such that
$$\{x\in K_z\ ;\ \langle x,z\rangle\ge0\}\subset \eps B_2^n, \quad \hbox{for all}\ |z|\ge M.$$
Moreover  $K_z$ contains the origin, hence
$$h_{K_{z}}\left(\frac{z}{|z|}\right)  = \max\left\{ \langle \frac{z}{|z|} ,v\rangle ;
v\in K_{z}, \langle z,v\rangle \ge 0\right\}\ra 0\ ,$$
when $|z|\ra +\infty$.
By part 2) of Lemma \ref{lemr},
$$h_{K_{z}}\left(-\frac{z}{|z|}\right)\ge r_{z}\left(-\frac{z}{|z|}\right)\ra +\infty\
.$$
\vskip 1mm\noindent
It follows that $\langle G(z),{\frac{z}{|z|}}\rangle \ra -\infty$, and
thus that $|G(z)|\ra +\infty$ when $|z|\ra +\infty$.
\vskip 1mm\noindent
But since $K_{z}$ is a convex body, $G(z)\in K_{z}$, and thus
$ |G(z)|\le r_{z}\left( \frac{G(z)}{ |G(z)|}\right).$
Since $|G(z)|\ra +\infty$,  one has $r_{z}\left(\frac{G(z)}{|G(z)|}\right)\ra +\infty$
when $|z|\ra +\infty$.
It follows again from part 1) of Lemma \ref{lemr} that for every
$\alpha<1$, there exists $M>0$ such that if $|z|>M$, then
$$\langle \frac{G(z)}{ |G(z)|},z \rangle \le -\alpha |z|\ . $$
This means that
$$\big\langle \frac{G(z)}{|G(z)|}, {\frac{z}{ |z|}}\big\rangle \ra 
-1\hbox{ when
}|z|\ra +\infty\ .$$

\vskip 1mm\noindent
{\bf B.} Let us prove that there exists $z_{0}\in \R^n$ such that $G(z_{0})=0$:
\vskip 1mm\noindent
Suppose that $G$ does not vanish.  Let
$C_{2}^n =\{x\in \R^n; |x|<1 \}$ be the open Euclidean unit ball, and define
$z: C_2^n\to\R^n$ by
$$
z(x):=\frac{x}{1-|x|} \ .
$$
Define also  $F: B_{2}^n\to S^{n-1}$ by
$$F(x)=\frac{G\left(z(x)\right)}{
|G\left(z(x)\right)|}\hbox{ for }x\in C_2^n,
\hbox{ and } F(u)=-u \hbox{ for  }u\in S^{n-1}.
$$
Let us prove that $F$ is continuous on $B_2^n$:
It is clear that
$F$ is continuous on $C_2^n$. Let  $u\in S^{n-1}$. 
If  $x\ra u$,  then $|z(x)|\ra +\infty$ and
$\frac{z(x)}{|z(x)|}=\frac{x }{|x|}\ra u$. Whence by {\bf A.},
$$ \big\langle \frac{z(x)}{ |z(x)|}, \frac{G\left(z(x)\right)}{
|G\left(z(x)\right)|}\big\rangle\ra -1,$$
which implies that 
$$F(x)=\frac{G\left(z(x)\right)}{ |G\left(z(x)\right)|}\ra -u\ .$$
Thus $F:B_{2}^n\to S^{n-1}$ is continuous and satisfies
$F(u)=-u$ for every $u\in S^{n-1}$. To conclude, we define
$Q:B_{2}^n\mapsto B_{2}^n$, by
$$Q(x)=\frac{x+F(x)}{
2}\hbox{ for every }x\in B_2^n\ .$$
Then $Q$ is continuous, but has
no fixed point, which  contradicts Brouwer fixed point theorem.
\hfill\qed

\vskip 3mm \noindent
{\bf Remark:}
\vskip 1mm \noindent
Theorem~\ref{Main} can be generalized in the following way: given 
$h:(0,+\infty)\to(0,+\infty)$ such that $t\mapsto h(e^t)$ is $\log$-concave 
and $h(r)r^{n-1}\to+\infty$ when $r\to+\infty$, let $\mu$ be the measure on $\R^n$ 
with density $h(|x|)$. 
Let $\rho:\R_{+}\to\R_{+}$ be measurable and
     $f:\R^n\to\R_+$ be
     a log-concave function such that $0<\int f d\mu<+\infty$. Then there
     exists $z\in \R^n$ such that for any measurable function $ g: \R^n \mapsto \R_{+}$
     satisfying
     $$f(x)g(y)\le\rho^2\left(\langle x-z, y-z\rangle\right)\ $$
for every  $x,y\in \R^n$ such that $ \langle x-z, y-z\rangle>0$,  one has
         $$\int_{\R^n} f(x)d\mu(x)\int_{\R^n} g(y)d\mu(y)\le
        \left( \int_{\R^n} \rho({|x|^{2}})d\mu(x)\right)^2.$$


   \section{Consequences on Legendre transform}

   Given a function $\phi:\R^n\to \R\cup\{+\infty\}$ and $z\in \R^n$, we
   recall that
the {\it Legendre transform}
$\L_{z}\phi$ of $\phi$ with
respect to $z\in \R^n$ is defined by
$$\L_{z}\phi(y) =\sup_{x}\left(\langle x-z,y-z\rangle-\phi(x)\right)\
\hbox{ for all } y\in \R^n\ .$$
For $z=0$, we use the notation $\L:=\L_0$. Observe that
$\L_{z}\phi:\R^n\to \R\cup\{+\infty\}$ is convex and that by a
classical separation argument, $\L_{z}(\L_{z}\phi)=\phi$, whenever
$\phi$ is itself convex
and $\phi(z)<+\infty$.  Notice also that the function $\phi(x)=|x|^{2}/2$ is
the unique function which satisfies $\L\phi=\phi$.
As a consequence of Theorem~\ref{Main}, we get the following theorem
which generalizes  the results of Artstein, Klartag and Milman \cite{AKM}
who  considered only the cases $\rho(t)=e^{-t}$ and $\rho(t)=(1-t)_+^m$.

\vskip 3mm\noindent
\begin{Thm}\label{Legendre} Let $\rho:\R_{+}\to\R_{+}$ be a log-concave
non-increasing function and let $\phi$ be a convex
function such that
$0<\int_{\R^n}\rho\left(\phi(x)\right)dx<+\infty$ . Then for some $z\in
\R^n$,
one has $$\int_{\R^n}\rho\left(\phi(x)\right)dx
\int_{\R^n}\rho\left(\L_{z}\phi(y)\right)dy\le
\left(\int_{\R^n}\rho\left(\frac{|x|^{2}}{2}\right)dx\right)^2.$$
If $\rho$ is decreasing, there is equality if and only if
for some positive definite matrix
$T:\R^n\to \R^n$ and
some $c\in \R$, one has
$$\phi(x)=\frac{|T(x+z)|^2}{ 2} +c\ ,\
\quad\hbox{ for all}\ x\in \R^n,$$
and moreover either $c=0$ or 
$\rho(t)=e^{at+b}$ for some $a< 0$, some $b\in\R$, and all
$t\in [-|c|, +\infty)$. 
\end{Thm}

\vskip 2mm\noindent
{\bf Proof.} 
      \vskip 1mm \noindent
{\bf A. The inequality.}
\vskip 1mm
\noindent
We apply Theorem~\ref{Main} to the $\log$-concave function
$f:=\rho\circ\phi$ to  get a convenient $z\in \R^n$.
By the definition of $\L_{z}$ and the fact that $\rho$ is $\log$-concave and
non-increasing, one has  for every $x,y\in \R^n$ such
that $\langle x-z,y-z\rangle>0$,
$$\rho\left(\phi(x)\right)\rho\left(\L_{z}\phi(y)\right)\le
\rho^2\left(\frac{\phi(x)+\L_{z}\phi(y)}{2}\right)
\le \rho^{2}\left(\frac{\langle x-z,y-z\rangle}{2}\right).$$
Setting $g(y)=\rho\left(\L_{z}\phi(y)\right)$, we may apply Theorem~\ref{Main}, to get the inequality.
\vskip 2mm\noindent
{\bf B. The case of equality.}
\vskip 1mm
\noindent
We may assume that
$z=0$. Set $\psi=\L\phi$. If there is equality,  we get from Theorem~\ref{Main} that for
some positive definite matrix
$T:\R^n\to \R^n$ and some $d>0$, one has
$$\frac{1}{ d}\rho\left( \phi(|T^{-1}x|)\right)=d
\rho\left(\psi(|Tx|)\right)=
\rho\left( \frac{|x|^2}{2} \right), $$
for every $x\in \R^n$. Since $\rho$ is $\log$-concave and
decreasing one has
\begin{eqnarray*}
\rho\left(\frac{|x|^{2}}{ 2}\right) & = &
\sqrt{\rho\left(\phi(T^{-1}x)\right)\rho\left(\psi(Tx)\right)}
\le \rho\left(\frac{\phi(T^{-1}x)+\psi(Tx)}{ 2}\right)\\
& \le & \rho\left(\frac{\langle T^{-1}x,Tx\rangle}{2}\right)=
\rho\left(\frac{|x|^{2}}{ 2}\right).
\end{eqnarray*}
Since $\rho$ is decreasing, we get
$\phi(T^{-1}x) +\psi(Tx)={|x|^2}\hbox{ for all $x\in \R^n$ .} $
Thus
$$
|x|^2-\phi(T^{-1}x)
 =  \psi(Tx)    =  \sup_{y}\big(\langle Tx,y\rangle-\phi(y)\big)
 =  \sup_{w}\left(\langle x,w\rangle-\phi(T^{-1}w)\right).
$$
We get $\phi(T^{-1}x)-\phi(T^{-1}w)\le |x|^{2} -\langle x,w\rangle$, for every $w,x\in \R^n$,
Setting $C(x)= \phi(T^{-1}x)- \frac{|x|^2}{ 2}$, it follows that
$$|C(x)-C(w)|\le \frac{|x-w|^{2}}{ 2}\hbox{ for all
$x,w\in \R^n$}. $$
It is easy then to conclude that $C$ is actually constant, and this
gives that for some $c>0$, one has
$$\phi(x)=\frac{|Tx|^2}{2} +c\quad {\rm and}\quad
\psi(x)=\frac{|T^{-1}x|^2}{ 2} -c \ .$$
This implies that $\rho$ satisfies
$$\rho\left(\frac{|x|^{2}}{ 2}\right)^2 =\rho\left(\frac{|x|^{2}}{ 2}+c\right)
\rho\left(\frac{|x|^{2}}{ 2}-c\right)
$$
and using again the log-concavity of $\rho$,
either $c=0$ or $\log(\rho)$ is affine on $[-|c|, +\infty)$.
\hfill\qed

\vskip 3mm\noindent
{\bf Remarks:}
\vskip 1mm\noindent
{\bf 1)} The cases when $\rho(t)=e^{-t}$ or $\rho(t)=(1-t)_+^m$ of
Theorem~\ref{Legendre}
were proved by Artstein, Klartag and Milman in \cite{AKM} by applying
the Blaschke-Santal\'o inequality for sets to a sequence of convex bodies
$(K_s(\phi))_{s\in\N}$
  in $\R^{n+s}$ and by letting $s\to +\infty$.  The use of this
  sequence makes the case of equality much more difficult that in our
  proof.

\vskip 1mm\noindent
{\bf 2)} In the case when the function $\rho$ is strictly convex 
(for example if $\rho(t)=e^{-t}$), then 
$$\min_z \int_{\R^n}\rho\left(\L_{z}\phi(y)\right)dy=
\min_z \int_{\R^n}\rho\left(\L\phi(y)-\langle z,y\rangle\right)dy$$
is reached at a unique point $z_{0}$ which satisfies 
$$z_{0}=\int_{\R^n}y\rho'(\L_{z_{0}}\phi(y))dy\bigg/
\int_{\R^n}\rho'(\L_{z_{0}}\phi(y))dy\ .$$
It follows that the inequality of Theorem~\ref{Legendre} is also
valid at this point $z=z_{0}$.

\vskip 1mm\noindent
{\bf 3)} Actually, it is also possible to prove Theorem \ref{Legendre}
by following step by step the method used by Meyer and Pajor (\cite{MP})
for proving Blaschke-Santal\'o inequality for convex bodies.
The idea is to prove that the quantity
$$\min_z \int_{\R^n}\rho\left(\L_{z}\phi(x)\right)dx$$
increases if we apply to  the epigraph
$E_\phi:=\{(x,t)\in\R^n\times\R\ ;\ \f(x)\le t\}$
of the function $\phi$ a well chosen Steiner symmetrisation to get a function
$\tilde{\phi}$ which is symmetric with respect ot the  symmetrisation hyperplane. 
After $n$ symmetrizations with respect to mutually orthogonal hyperplanes,
the function is unconditional and the result follows from
the application of the Pr\'ekopa-Leindler inequality for the geometric mean
(Theorem \ref{Prekopageom}).
However, this proof is much longer, and seems to require some  additionally
hypotheses on the function $\rho$, namely that $\rho$ is convex and decreasing and that $-\rho'$ is
$\log$-concave.

\vskip 1mm\noindent
{\bf  4)  Shortcut for the proof of the equality case in
Blaschke-Santal\'o inequality.}

There exists different proofs of the equality case for
Blaschke-Santalo's inequality. It was first proved in the centrally
symmetric case by Saint-Raymond \cite{Saint-Raymond},
using a tricky lemma for functions of one variable,
then in the general case by Petty
\cite{Petty} with some involved arguments of PDE (see also D.~Hug \cite{Hug}). A simpler proof together with a stronger
inequality was then given by Meyer and Pajor  \cite{MP} using the Steiner symmetrization,
a result of \cite{Falconer} and finally the lemma of Saint-Raymond.
\vskip 1mm\noindent
In fact, one can give  the following simpler argument.
\vskip 1mm\noindent
{\bf a.} If
$K$ is unconditional with maximal volume product, we have seen that the case of
equality follows easily from the equality case
in the one-dimensional Pr\'ekopa-Leindler inequality.
\vskip 1mm\noindent
{\bf b.} Suppose now that $K$ has maximal volume product and is centrally symmetric.
Then for every  $u\in S^{n-1}$, after $n$  Steiner symmetrizations 
with respect to pairwise orthogonal hyperplanes, the last one being
with respect to $\{u\}^{\perp}$, we get from $K$ an unconditional body
with maximal volume product (recall that a Steiner symmetrization does
not decrease volume product), and thus by {\bf a.} an ellipsoid. To 
conclude that
$K$ is itself an ellipsoid, we use the following elementary lemma, where for
$v\in S^{n-1}$, we denote by $S_{v}K$ the Steiner symmetral of $K$
with respect to the hyperplane $v^{\perp}:=\{x\in\R^n\ ;\ \langle x, v\rangle=0\}$.
\vskip 3mm\noindent
{\bf Lemma.} Let $K$ be a centrally symmetric convex body. 
Then $K$ is an ellipsoid if and only if  for
every orthonormal basis $(u_{1},\dots,u_{n})$ of $\R^n$,
$S_{u_n}S_{u_{n-1}}\dots S_{u_1}K$ is an ellipsoid.
\vskip 3mm\noindent
{\bf Proof:} The "only if " part is well known. For the "if" part, fix $u\in  S^{n-1}$, 
and $(u_{1},\dots,u_{n})$ be an orthonormal basis such that $u=u_{n}$. 
Let $L=S_{u_{n-1}}\cdots S_{u_1}K$. Then $L$
is centrally symmetric (since $K$ is), and symmetric
with respect to the $(n-1)$ pairwise orthogonal hyperplanes
$u_{i}^{\perp}$, $1\le i\le n-1$. It follows that $L$ is also
symmetric with respect to  $u_{n}^{\perp}$, so that
$L= S_{u_{n}}L= S_{u_n}S_{u_{n-1}}\cdots S_{u_1}K$ is an ellipsoid.
Thus for some $a_1,\dots, a_n>0$ one has
$$L=\left\{x=x_{1}u_{1}+\cdots +x_{n}u_{n}; \sum_{i=1}^n
\frac{x_{n}^{2}}
{a_{n}^{2}}\le 1\right\}.$$
Let $h_{K}(u):= \max\{\langle x,u\rangle ;  x\in K\} $.
It is easy to see that whenever $v\in S^{n-1}$ satisfy
$\langle v,u\rangle =0$, then 
$$h_{K}(u)=h_{S_{v}K}(u)\ \hbox{  and }\ 
\int_{K}\langle x,u\rangle^{2}dx= \int_{S_{v}K}\langle x,u\rangle^{2}dx.$$
It follows that $a_{n}=h_{L}(u_{n})=h_{K}(u_{n})$ and  
$$\int_{K}\langle x,u_{n}\rangle^{2}dx=
\int_{L}\langle x,u_{n}\rangle^{2}dx =\frac{v_{n}}{ n+2}\cdot
a_{1}\dots a_{n}\cdot a_{n}^2\ .$$
Since $|L|=|K|$, one has $v_{n}a_{1}\cdots a_{n}=|K|$. Thus
$$h_{K}(u)^{2}=\frac{n+2}{ |K|} \int_{K}\langle
x,u\rangle^{2}dx\hbox{ for every } u\in S^{n-1}. $$
It follows that $K^\circ$ and thus $K$ is an ellipsoid.  
\hfill \qed

\bigskip


%

\end{document}